\documentclass[12pt]{article}
\usepackage{latexsym,amsfonts,amssymb}
\usepackage[dvips]{color}
\usepackage{colordvi,multicol}

\def \cal{\mathcal}

\newtheorem{thm}{Theorem}[section]
\newtheorem{cor}[thm]{Corollary}

\newtheorem{rem}[thm]{Remark}

\begin{document}
\title{\bf  A system of nonlinear equations with application to large deviations for Markov chains with finite lifetime}
\author{}

\maketitle
\date{}

 \centerline{Ze-Chun Hu} \centerline{\small College
of Mathematics, Sichuan University, Chengdu, 610064, China}
\centerline{\small E-mail: zchu@scu.edu.cn}
\vskip 0.7cm \centerline{Wei Sun} \centerline{\small Department of
Mathematics and Statistics, Concordia University,}
\centerline{\small Montreal, H3G 1M8, Canada} \centerline{\small
E-mail: wei.sun@concordia.ca}
\vskip 0.7cm \centerline{Jing Zhang} \centerline{\small School of
Mathematics and Statistics, Hainan Normal University,}
\centerline{\small  Haikou, 571158, China}
\centerline{\small
E-mail: zh\_jing0820@hotmail.com}

\vskip 1cm


\vskip 0.5cm \noindent{\bf Abstract}\quad In this paper, we first show the existence of solutions to the following system of nonlinear equations
\begin{eqnarray*}
\left\{
\begin{array}{l}
a_{11}x_1+a_{12}x_2+a_{13}x_3+\cdots+a_{1n}x_{n}
=b_{11}\frac{1}{x_1}+b_{12}\frac{1}{x_2}+b_{13}\frac{1}{x_3}+\cdots+b_{1n}\frac{1}{x_{n}},\\
a_{21}\frac{1}{x_1}+a_{22}\frac{x_2}{x_1}+a_{23}\frac{x_3}{x_1}+\cdots+a_{2n}\frac{x_{n}}{x_1}
=b_{21}x_1+b_{22}\frac{x_1}{x_2}+b_{23}\frac{x_1}{x_3}+\cdots+b_{2n}\frac{x_1}{x_{n}},\\
a_{31}\frac{x_1}{x_2}+a_{32}\frac{1}{x_2}+a_{33}\frac{x_3}{x_2}+\cdots+a_{3n}\frac{x_{n}}{x_2}
=b_{31}\frac{x_2}{x_1}+b_{32}x_2+b_{33}\frac{x_2}{x_3}+\cdots+b_{3n}\frac{x_2}{x_{n}},\\
\cdots\cdots\\
a_{n1}\frac{x_1}{x_{n-1}}+a_{n2}\frac{x_2}{x_{n-1}}+a_{n3}\frac{x_3}{x_{n-1}}+
\cdots+a_{n,n-1}\frac{1}{x_{n-1}}+a_{nn}\frac{x_{n}}{x_{n-1}}\\
=b_{n1}\frac{x_{n-1}}{x_1}+b_{n2}\frac{x_{n-1}}{x_2}+b_{n3}\frac{x_{n-1}}{x_3}+\cdots+b_{n,n-1}x_{n-1}
+b_{nn}\frac{x_{n-1}}{x_{n}},
\end{array}
\right.
\end{eqnarray*}
where $n\ge 3$ and $a_{ij},b_{ij},1\le i,j\le n$,
are positive constants. Then, we make use of this result to obtain the large deviation principle
for the occupation time distributions of continuous-time finite state Markov chains with finite lifetime.

\smallskip


\section{Introduction and main results}

In a series of fundamental papers (see \cite{I,II,III,IV}), Donsker and Varadhan developed the large deviation theory for the occupation time distributions of Markov processes. By virtue of Dirichlet forms, Fukushima and Takeda derived the Donsker-Varadhan type large deviation principle for a general, not necessarily conservative symmetric Markov processes (see \cite{fukutakeda}, \cite[Section 6.4]{fukushima2010} and the references therein). The motivation of this work is to generalize some results of Donsker-Varadhan and Fukushima-Takeda to not necessarily conservative and not necessarily symmetric Markov processes.

We denote $E=\{1,2,\dots, n\}$ for $n\ge 2$. Let $X=((X_t)_{t\ge 0},(P_i)_{i\in E})$ be a continuous-time Markov chain with the state space $E$. Denote by
$\zeta$ the lifetime of $X$ and denote by $Q=(q_{ij})_{i,j\in E}$ the $Q$-matrix of $X$. We assume that $Q$ satisfies the following conditions:

(1) $0<-q_{ii}<\infty,\ i\in E$.

(2) $q_{ij}>0,\ i,j\in E,\ i\neq j$.

(3) $\sum_{j\in E}q_{ij}\le0,\ i\in E$.

\noindent In this paper, we will derive the large deviation principle for the occupation time distributions of $X$.

We discover that the large deviations for $X$ rely heavily on the existence of solutions to the following system of nonlinear equations
\begin{eqnarray}\label{sys}
\left\{
\begin{array}{l}
a_{11}x_1+a_{12}x_2+a_{13}x_3+\cdots+a_{1n}x_{n}
=b_{11}\frac{1}{x_1}+b_{12}\frac{1}{x_2}+b_{13}\frac{1}{x_3}+\cdots+b_{1n}\frac{1}{x_{n}},\\
a_{21}\frac{1}{x_1}+a_{22}\frac{x_2}{x_1}+a_{23}\frac{x_3}{x_1}+\cdots+a_{2n}\frac{x_{n}}{x_1}
=b_{21}x_1+b_{22}\frac{x_1}{x_2}+b_{23}\frac{x_1}{x_3}+\cdots+b_{2n}\frac{x_1}{x_{n}},\\
a_{31}\frac{x_1}{x_2}+a_{32}\frac{1}{x_2}+a_{33}\frac{x_3}{x_2}+\cdots+a_{3n}\frac{x_{n}}{x_2}
=b_{31}\frac{x_2}{x_1}+b_{32}x_2+b_{33}\frac{x_2}{x_3}+\cdots+b_{3n}\frac{x_2}{x_{n}},\\
\cdots\cdots\\
a_{n1}\frac{x_1}{x_{n-1}}+a_{n2}\frac{x_2}{x_{n-1}}+a_{n3}\frac{x_3}{x_{n-1}}+
\cdots+a_{n,n-1}\frac{1}{x_{n-1}}+a_{nn}\frac{x_{n}}{x_{n-1}}\\
=b_{n1}\frac{x_{n-1}}{x_1}+b_{n2}\frac{x_{n-1}}{x_2}+b_{n3}\frac{x_{n-1}}{x_3}+\cdots+b_{n,n-1}x_{n-1}
+b_{nn}\frac{x_{n-1}}{x_{n}},
\end{array}
\right.
\end{eqnarray}
where $n\ge 3$ and $a_{ij},b_{ij},1\le i,j\le n$,
are constants. It is a bit surprising to us that (\ref{sys}) turns out to be undiscussed to date. In the next section, we will prove the following result.
\begin{thm}\label{thm1}
Suppose that $n\ge 3$ and $a_{ij},b_{ij},1\le i,j\le n$,
are positive constants. Then, there exists a positive solution $(x_1,x_2,\dots,x_n)$ to (\ref{sys}).
\end{thm}

As a direct consequence of Theorem \ref{thm1}, we obtain the following result.
\begin{thm}\label{thm2} Suppose that $n\ge 2$ and $\beta(i)>0, 1\le i\le n$. Then, there exist $\alpha(i)>0, 1\le i\le n$, such that
\begin{equation}\label{hsz1}
\sum_{j\in E}q_{ij}\frac{\alpha(j)}{\alpha(i)}\beta(j)=\sum_{j\in
E}q_{ji}\frac{\alpha(i)}{\alpha(j)}\beta(j),\ \ \forall i\in E.
\end{equation}
\end{thm}
The proof of Theorem \ref{thm2} will be given in the next section.
\begin{rem}
(a) Denote by $A$ the
diagonal matrix with $A_{ii}=\alpha(i)$, $1\le i \le n$,
$A_{ij}=0$ if $i\not=j$, and denote by
$\beta=(\beta(1),\dots,\beta(n))^T$. Hereafter $^T$
denotes transpose. Then, we can rewrite (\ref{hsz1}) as
follows
\begin{equation}\label{rem33}
A^{-1}QA\beta=AQ^TA^{-1}\beta.
\end{equation}
Theorem \ref{thm2} implies that for
any vector $\beta>0$ there exists a positive
diagonal matrix $A$ such that (\ref{rem33}) holds.

\noindent (b) If the
matrix $Q$ is symmetric, then it is easy to see that
$\alpha(i)\equiv1, 1\le i\le n$, provide a solution to (\ref{hsz1}). When $Q$ is non-symmetric, Theorem \ref{thm2} seems to be a new result in the literature.
\end{rem}

In Section 3 of this paper, we will make use of Theorem \ref{thm2} to obtain the
large deviation principle for $X$. Define the normalized occupation time distribution
$L_t$, $t>0$, by
$$
L_t(A)=\frac{1}{t}\int_0^t 1_A(X_s)ds,\ \ A\subset E.
$$
Let $u$ be a function  on $E$. We write $u=(u(1),\dots,u(n))^T$
and denote $u>0$ if $u(i)>0$ for each $i\in E$. For $i\in E$, we have $Qu(i)=\sum_{j\in E}q_{ij}u(j)$. Let
$\mu$ be a measure on $E$. We define
$$
I(\mu):=-\inf_{u>0}\int_E\frac{Qu}{u}d\mu=-\inf_{u>0}\sum_{i\in
E}\frac{Qu(i)}{u(i)}\mu(\{i\}).
$$
Denote by ${\cal P}_1(E)$ the set of all probability measures on
$E$.

\begin{thm}\label{hsz3} For each open set $G$ of ${\cal
P}_1(E)$, \begin{equation}\label{new1}
\liminf_{t\rightarrow\infty}\frac{1}{t}\log P_{i}(L_t\in G,\
t<\zeta)\ge -\inf_{\mu\in G}I(\mu),\ \ \forall i\in E.
\end{equation}
For each closed set $K$ of ${\cal P}_1(E)$, \begin{equation}\label{lower}
\limsup_{t\rightarrow\infty}\frac{1}{t}\log P_{i}(L_t\in K,\
t<\zeta)\le -\inf_{\mu\in K}I(\mu),\ \ \forall i\in E.
\end{equation}
\end{thm}

\vskip 0.2cm By setting $G=K={\cal P}_1(E)$ in Theorem \ref{hsz3},
we get
\begin{cor} For $i\in E$,
$$\lim_{t\rightarrow\infty}\frac{1}{t}\log P_{i}(t<\zeta)=\sup_{\mu\in {\cal P}_1(E)}\left\{\inf_{u>0}\sum_{i\in
E}\sum_{j\in E}\frac{q_{ij}\mu(\{i\})u(j)}{u(i)}\right\}.
$$
\end{cor}

\section{Proofs of Theorems \ref{thm1} and \ref{thm2}}\setcounter{equation}{0}
\noindent {\bf Proof of Theorem \ref{thm1}}.
\vskip 0.2cm

\emph{We first consider the case that $n=3$.}
\vskip 0.2cm

We define four continuous functions $f_1,f_2,f_3,F$ with the
domain $D_3:=\{x=(x_1,x_2,x_3)\in \mathbf{R}^3|x_i>0,1\le i\le 3\}$ by
\begin{eqnarray*}
&&f_1(x)=(a_{11}x_1+a_{12}x_2+a_{13}x_3)-\left(b_{11}\frac{1}{x_1}+b_{12}\frac{1}{x_2}+b_{13}\frac{1}{x_3}\right),\\
&&f_2(x)=\left(a_{21}\frac{1}{x_1}+a_{22}\frac{x_2}{x_1}+a_{23}\frac{x_3}{x_1}\right)
-\left(b_{21}x_1+b_{22}\frac{x_1}{x_2}+b_{23}\frac{x_1}{x_3}\right),\\
&&f_3(x)=\left(a_{31}\frac{x_1}{x_2}+a_{32}\frac{1}{x_2}+a_{33}\frac{x_3}{x_2}\right)
-\left(b_{31}\frac{x_2}{x_1}+b_{32}x_2+b_{33}\frac{x_2}{x_3}\right),\\
&&F(x)=f_1^2(x)+f_2^2(x)+f_3^2(x).
\end{eqnarray*}
It is easy to see that the function $F$ has a minimum value at some point $x^*=(x_1^*,x_2^*,x_3^*)\in D_3$.  In the following, we will prove that
$$
f_1(x^*)=f_2(x^*)=f_3(x^*)=0.
$$

Since $F$ has a minimum value at $x^*$, we have
\begin{eqnarray}
&&\frac{\partial F}{\partial x_1}(x^*)\!=\!f_1(x^*)c_{11}\!-\!f_2(x^*)c_{12}\!+\!f_3(x^*)c_{13}\!=\!0,\label{4}\\
&&\frac{\partial F}{\partial x_2}(x^*)\!=\!f_1(x^*)c_{21}\!+\!f_2(x^*)c_{22}\!-\!f_3(x^*)c_{23}\!=\!0,\label{5}\\
&&\frac{\partial F}{\partial x_3}(x^*)\!=\!f_1(x^*)c_{31}\!+\!f_2(x^*)c_{32}\!+\!f_3(x^*)c_{33}\!=\!0,\label{6}
\end{eqnarray}
where $c_{ij},1\le i,j\le 3$, are positive constants. If
$f_3(x^*)=0$, then we obtain by (\ref{4}) and (\ref{5})
that $f_1(x^*)=f_2(x^*)=0$. Similarly, if
$f_1(x^*)=0$ or $f_2(x^*)=0$, we also have
$$
f_1(x^*)=f_2(x^*)=f_3(x^*)=0.
$$
Thus, to prove the existence of solutions to (\ref{sys}), we need only show that there is a contradiction if
\begin{equation}\label{r1}
f_1(x^*)*f_2(x^*)*f_3(x^*)\neq 0.
\end{equation}

Suppose that (\ref{r1}) holds. If $f_1(x^*)>0$, then we obtain by (\ref{6}) that either $f_2(x^*)<0$ or $f_3(x^*)<0$. Further, we obtain by (\ref{4}) and (\ref{5}) that both $f_2(x^*)<0$ and $f_3(x^*)<0$. Similarly, we can show that if $f_1(x^*)<0$, then $f_2(x^*)>0$ and $f_3(x^*)>0$. Therefore, to prove the existence of solutions to (\ref{sys}), we need only show that neither of the following two cases can happen:
\vskip 0.2cm
{\bf Case (i).} $f_1(x^*)>0$, $f_2(x^*)<0$ and $f_3(x^*)<0$.
\vskip 0.1cm
{\bf Case (ii). } $f_1(x^*)<0$, $f_2(x^*)>0$ and $f_3(x^*)>0$.
\vskip.2cm
{\bf Case (i) cannot happen.}  Suppose that
\begin{eqnarray*}
&&f_1(x^*)=(a_{11}x^*_1+a_{12}x^*_2+a_{13}x^*_3)-\left(b_{11}\frac{1}{x^*_1}+b_{12}\frac{1}{x^*_2}+b_{13}\frac{1}{x^*_3}\right)>0,\\
&&f_2(x^*)=\left(a_{21}\frac{1}{x^*_1}+a_{22}\frac{x^*_2}{x^*_1}+a_{23}\frac{x^*_3}{x^*_1}\right)
-\left(b_{21}x^*_1+b_{22}\frac{x^*_1}{x^*_2}+b_{23}\frac{x^*_1}{x^*_3}\right)<0,\\
&&f_3(x^*)=\left(a_{31}\frac{x^*_1}{x^*_2}+a_{32}\frac{1}{x^*_2}+a_{33}\frac{x^*_3}{x^*_2}\right)
-\left(b_{31}\frac{x^*_2}{x^*_1}+b_{32}x^*_2+b_{33}\frac{x^*_2}{x^*_3}\right)<0.
\end{eqnarray*}
In the following, we will show that there exist sufficiently small positive numbers $\delta_1$ and $\delta_2$ such that
$x^*_1-\delta_1>0,x^*_2-\delta_2>0$, $
\frac{x^*_1}{x^*_2}=\frac{\delta_1}{\delta_2}
$, and
\begin{eqnarray*}
&&f_1(x^*)>f_1(x^*_1-\delta_1,x^*_2-\delta_2,x^*_3)>0,\\
&&f_2(x^*)<f_2(x^*_1-\delta_1,x^*_2-\delta_2,x^*_3)<0,\\
&&f_3(x^*)<f_3(x^*_1-\delta_1,x^*_2-\delta_2,x^*_3)<0.
\end{eqnarray*}

Define $\delta=\frac{x^*_2}{x^*_1},\delta_2=\delta\delta_1$, and
$$
\gamma_1=f_1(x^*),\ \gamma_2=-f_2(x^*),\  \gamma_3=-f_3(x^*).
$$
Then, it is sufficient to show that there exists a positive number $\delta_1$ such that
\begin{eqnarray*}
&&\delta_1<x^*_1,\ \delta\delta_1<x^*_2,\\
&&0<f_1(x^*)-f_1(x^*_1-\delta_1,x^*_2-\delta_2,x^*_3)\leq \frac{\gamma_1}{2},\\
&&0<f_2(x^*_1-\delta_1,x^*_2-\delta_2,x^*_3)-f_2(x^*)\leq \frac{\gamma_2}{2},\\
&&0<f_3(x^*_1-\delta_1,x^*_2-\delta_2,x^*_3)-f_3(x^*)\leq \frac{\gamma_3}{2},
\end{eqnarray*}
i.e.,
\begin{eqnarray*}
&&\delta_1<\min\left\{x^*_1,\frac{x^*_2}{\delta}\right\},\\
&&0<(a_{11}+a_{12}\delta)\delta_1+\left(b_{11}+\frac{b_{12}}{\delta}\right)\left(\frac{1}{x^*_1-\delta_1}
-\frac{1}{x^*_1}\right)\leq \frac{\gamma_1}{2},\\
&&0<(a_{21}+a_{23}x^*_3)\left(\frac{1}{x^*_1-\delta_1}-\frac{1}{x^*_1}\right)+\left(b_{21}+\frac{b_{23}}{x^*_3}\right)
\delta_1\leq \frac{\gamma_2}{2},\\
&&0<(a_{32}+a_{33}x^*_3)\left(\frac{1}{\delta}\right)\left(\frac{1}{x^*_1-\delta_1}-\frac{1}{x^*_1}\right)
+\left(b_{32}+\frac{b_{33}}{x^*_3}\right)\delta\delta_1\leq \frac{\gamma_3}{2}.
\end{eqnarray*}
Obviously, there exists a positive number $\delta_1$ satisfying all the above conditions. For this $\delta_1$, we have that $F(x^*_1-\delta_1,x^*_2-\delta_2,x^*_3)<F(x^*)$, which contradicts that $F$ reaches its minimum at $x^*$.
\vskip 0.2cm
{\bf Case (ii) cannot happen.}  Suppose that
\begin{eqnarray*}
&&f_1(x^*)=(a_{11}x^*_1+a_{12}x^*_2+a_{13}x^*_3)-\left(b_{11}\frac{1}{x^*_1}+b_{12}\frac{1}{x^*_2}+b_{13}\frac{1}{x^*_3}\right)<0,\\
&&f_2(x^*)=\left(a_{21}\frac{1}{x^*_1}+a_{22}\frac{x^*_2}{x^*_1}+a_{23}\frac{x^*_3}{x^*_1}\right)
-\left(b_{21}x^*_1+b_{22}\frac{x^*_1}{x^*_2}+b_{23}\frac{x^*_1}{x^*_3}\right)>0,\\
&&f_3(x^*)=\left(a_{31}\frac{x^*_1}{x^*_2}+a_{32}\frac{1}{x^*_2}+a_{33}\frac{x^*_3}{x^*_2}\right)
-\left(b_{31}\frac{x^*_2}{x^*_1}+b_{32}x^*_2+b_{33}\frac{x^*_2}{x^*_3}\right)>0.
\end{eqnarray*}
In the following, we will show that there exist sufficiently small positive numbers $\delta_1$ and $\delta_2$ such that $
\frac{x^*_1}{x^*_2}=\frac{\delta_1}{\delta_2}$,
and
\begin{eqnarray*}
&&f_1(x^*)<f_1(x^*_1+\delta_1,x^*_2+\delta_2,x^*_3)<0,\\
&&f_2(x^*)>f_2(x^*_1+\delta_1,x^*_2+\delta_2,x^*_3)>0,\\
&&f_3(x^*)>f_3(x^*_1+\delta_1,x^*_2+\delta_2,x^*_3)>0.
\end{eqnarray*}

Define $\delta=\frac{x^*_2}{x^*_1},\delta_2=\delta\delta_1$, and
$$
\gamma_1=-f_1(x^*),\ \gamma_2=f_2(x^*),\  \gamma_3=f_3(x^*).
$$
Then, it is sufficient to show that there exists a positive number $\delta_1$ such that
\begin{eqnarray*}
&&0<f_1(x^*_1+\delta_1,x^*_2+\delta_2,x^*_3)-f_1(x^*)\leq \frac{\gamma_1}{2},\\
&&0<f_2(x^*)-f_2(x^*_1+\delta_1,x^*_2+\delta_2,x^*_3)\leq \frac{\gamma_2}{2},\\
&&0<f_3(x^*)-f_3(x^*_1+\delta_1,x^*_2+\delta_2,x^*_3)\leq \frac{\gamma_3}{2},
\end{eqnarray*}
i.e.,
\begin{eqnarray*}
&&0<(a_{11}+a_{12}\delta)\delta_1+\left(b_{11}+\frac{b_{12}}{\delta}\right)\left(\frac{1}{x^*_1}
-\frac{1}{x^*_1+\delta_1}\right)\leq \frac{\gamma_1}{2},\\
&&0<\left(b_{21}+\frac{b_{23}}{x^*_3}\right)\delta_1+(a_{21}+a_{23}x^*_3)\left(\frac{1}{x^*_1}-\frac{1}{x^*_1+\delta_1}\right)\leq \frac{\gamma_2}{2},\\
&&0<\left(b_{32}+\frac{b_{33}}{x^*_3}\right)\delta\delta_1+(a_{32}+a_{33}x^*_3)\left(\frac{1}{\delta}\right)
\left(\frac{1}{x^*_1}-\frac{1}{x^*_1+\delta_1}\right)\leq \frac{\gamma_3}{2}.
\end{eqnarray*}
Obviously, there exists a positive number $\delta_1$ satisfying
all the above conditions. For this $\delta_1$, we have  that
$F(x^*_1+\delta_1,x^*_2+\delta_2,x^*_3)<F(x^*)$, which
contradicts that $F$ reaches its minimum at $x^*$.
\vskip 0.2cm

\emph{We now consider the general case that $n\ge 4$.}
\vskip 0.2cm

We define $(n+1)$ continuous functions $f_1,f_2,\dots,f_{n},F$ with
the domain $D_{n}:=\{x=(x_1,x_2,\dots,x_{n}) \in
\mathbf{R}^{n}|x_i>0,i=1,2,\dots,n\}$  by
\begin{eqnarray*}
&&f_1(x)=\left(a_{11}x_1+a_{12}x_2+\cdots+a_{1n}x_{n}\right)\\
&&\hspace{4cm}-\left(b_{11}\frac{1}{x_1}+b_{12}\frac{1}{x_2}+\cdots+b_{1n}\frac{1}{x_{n}}\right),\\
&&f_2(x)=\left(a_{21}\frac{1}{x_1}+a_{22}\frac{x_2}{x_1}+\cdots+a_{2n}\frac{x_{n}}{x_1}\right)\\
&&\hspace{4cm}-\left(b_{21}x_1+b_{22}\frac{x_1}{x_2}+\cdots+b_{2n}\frac{x_1}{x_{n}}\right),\\
&&\cdots\cdots\\
&&f_{n}(x)=\left(a_{n,1}\frac{x_1}{x_{n-1}}+a_{n,2}\frac{x_2}{x_{n-1}}+\cdots+a_{n,n}\frac{x_{n}}{x_{n-1}}\right)\\
&&\hspace{4cm}-\left(b_{n,1}\frac{x_{n-1}}{x_1}+b_{n,2}\frac{x_{n-1}}{x_2}+\cdots+b_{n,n}\frac{x_{n-1}}{x_{n}}\right),\\
&&F(x)=\sum_{i=1}^{n}f_i^2(x).
\end{eqnarray*}
The function $F$ has a minimum value at some
point $x^*=(x_1^*,x_2^*,\dots,x_{n}^*)\in D_{n}$.  In the
following, we will prove that
$$
f_i(x^*)= 0,\ \ \forall i\in \{1,2,\dots,n\}.
$$

Since $F$ has a minimum value at $x^*$, we have
$$
\frac{\partial F}{\partial x_i}(x^*)=0,\ \ \forall i\in
\{1,2,\dots,n\}.
$$
It follows that
\begin{eqnarray*}
\left\{
\begin{array}{l}
f_1(x^*)c_{11}-f_2(x^*)c_{12}+f_3(x^*)c_{13}+\cdots+f_{n}(x^*)c_{1n}=0,\\
f_1(x^*)c_{21}+f_2(x^*)c_{22}-f_3(x^*)c_{23}+\cdots+f_{n}(x^*)c_{2n}=0,\\
\cdots\cdots\\
f_1(x^*)c_{n-1,1}+\cdots+f_{n-1}(x^*)c_{n-1,n-1}-f_{n}(x^*)c_{n-1,n}=0,\\
f_1(x^*)c_{n1}+f_2(x^*)c_{n2}+\cdots+f_{n}(x^*)c_{nn}=0,
\end{array}
\right.
\end{eqnarray*}
where $c_{ij},1\le i,j\le n$,
are positive constants. Note that there is exactly one minus sign
in the first $(n-1)$ equations and there is no minus sign in the
last equation.
\vskip 0.2cm
{\bf Case (a).} Suppose that $f_1(x^*)=0$. We
consider the following $(n+1)$ continuous functions with the domain
$\bar{D}_{n-1}=\{(x_2,\dots,x_{n})\in
\mathbf{R}^{n-1}|x_i>0,i=2,\dots,n\}$:
\begin{eqnarray*}
&&\bar{f}_i(x_2,\dots,x_{n}):=f_i(x_1^*,x_2,\dots,x_{n}),\ \ i=1,\dots,n,\\
&&\bar{F}(x_2,\dots,x_{n}):=F(x_1^*,x_2,\dots,x_{n}).
\end{eqnarray*}
Since $F$ has a minimum value at $x^*\in
D_{n}$, $\bar{F}$ has a minimum value at
$(x_2^*,\dots,x_{n}^*)\in\bar{D}_{n-1}$.  Then, we have
$$
\frac{\partial \bar{F}}{\partial x_i}(x_2^*,\dots,x_{n}^*)=0,\ \
i=2,\dots,n,
$$
which together with $\bar{f}_1(x_2^*,\dots,x_{n}^*)=f_1(x^*)=0$
implies that
\begin{eqnarray*}
\left\{
\begin{array}{l}
\bar{f}_2(x_2^*,\dots,x_{n}^*)\bar{c}_{22}-\bar{f}_3(x_2^*,\dots,x_{n}^*)\bar{c}_{23}
+\bar{f}_4(x_2^*,\dots,x_{n}^*)\bar{c}_{24}+\cdots+\bar{f}_{n}(x_2^*,\dots,x_{n}^*)\bar{c}_{2n}=0,\\
\bar{f}_2(x_2^*,\dots,x_{n}^*)\bar{c}_{32}+\bar{f}_3(x_2^*,\dots,x_{n}^*)\bar{c}_{33}
-\bar{f}_4(x_2^*,\dots,x_{n}^*)\bar{c}_{34}+\cdots+\bar{f}_{n}(x_2^*,\dots,x_{n}^*)\bar{c}_{3,n}=0,\\
\cdots\cdots\\
\bar{f}_2(x_2^*,\dots,x_{n}^*)\bar{c}_{n-1,2}+\cdots+\bar{f}_{n-1}(x_2^*,\dots,x_{n}^*)\bar{c}_{n-1,n-1}
-\bar{f}_{n}(x_2^*,\dots,x_{n}^*)\bar{c}_{n-1,n}=0,\\
\bar{f}_2(x_2^*,\dots,x_{n}^*)\bar{c}_{n2}+\bar{f}_3(x_2^*,\dots,x_{n}^*)\bar{c}_{n3}
+\cdots+\bar{f}_{n}(x_2^*,\dots,x_{n}^*)\bar{c}_{nn}=0,
\end{array}
\right.
\end{eqnarray*}
where $\bar{c}_{ij},1\le i,j\le n$, are positive
constants. Thus, we obtain by following the same argument for the
$(n-1)$ case that
$$
\bar{f}_i(x_2^*,\dots,x_{n}^*)=0,\ \ i=2,\dots,n.
$$
Therefore,
$$
f_i(x^*)= 0,\ \ \forall i\in \{1,2,\dots,n\}.
$$
\vskip 0.2cm
{\bf Case (b).} Suppose that
$\prod_{i=2}^{n}f_i(x^*)=0$. By
symmetry, we can assume without loss of generality that
$f_{n}(x^*)=0$. Now we consider the
following $(n+1)$ continuous functions with the domain
$\bar{D}_{n-1}=\{(x_1,\dots,x_{n-1})\in
\mathbf{R}^{n-1}|x_i>0,i=1,\dots,n-1\}$:
\begin{eqnarray*}
&&\bar{f}_i(x_1,\dots,x_{n-1}):=f_i(x_1,\dots,x_{n-1},x^*_{n}),\ \ i=1,\dots,n,\\
&&\bar{F}(x_1,\dots,x_{n-1}):=F(x_1,\dots,x_{n-1},x_{n}^*).
\end{eqnarray*}
Since $F$ has a minimum value at $x^*\in
D_{n}$, $\bar{F}(x_1,\dots,x_{n-1})$ has a minimum value at
$(x_1^*,\dots,x_{n-1}^*)\in\bar{D}_{n-1}$.  Then, we have
$$
\frac{\partial \bar{F}}{\partial x_i}(x_1^*,\dots,x_{n-1}^*)=0,\ \
i=1,\dots,n-1,
$$
which together with
$\bar{f}_{n}(x_1^*,\dots,x_{n-1}^*)=f_{n}(x^*)=0$ implies that
\begin{eqnarray*}
\left\{
\begin{array}{l}
\bar{f}_1(x_1^*,\dots,x_{n-1}^*)\bar{c}_{11}-\bar{f}_2(x_1^*,\dots,x_{n-1}^*)\bar{c}_{12}
+\bar{f}_3(x_1^*,\dots,x_{n-1}^*)\bar{c}_{13}+\cdots+\bar{f}_{n-1}(x_1^*,\dots,x_{n-1}^*)\bar{c}_{1,n-1}=0,\\
\bar{f}_1(x_1^*,\dots,x_{n-1}^*)\bar{c}_{21}+\bar{f}_2(x_1^*,\dots,x_{n-1}^*)\bar{c}_{22}
-\bar{f}_3(x_1^*,\dots,x_{n-1}^*)\bar{c}_{23}+\cdots+\bar{f}_{n-1}(x_1^*,\dots,x_{n-1}^*)\bar{c}_{2,n-1}=0,\\
\cdots\cdots\\
\bar{f}_1(x_1^*,\dots,x_{n-1}^*)\bar{c}_{n-2,1}+\cdots+\bar{f}_{n-2}(x_1^*,\dots,x_{n-1}^*)\bar{c}_{n-2,n-2}
-\bar{f}_{n-1}(x_1^*,\dots,x_{n-1}^*)\bar{c}_{n-2,n-1}=0,\\
\bar{f}_1(x_1^*,\dots,x_{n-1}^*)\bar{c}_{n-1,1}+\bar{f}_2(x_1^*,\dots,x_{n-1}^*)\bar{c}_{n-1,2}+\cdots
+\bar{f}_{n-1}(x_1^*,\dots,x_{n-1}^*)\bar{c}_{n-1,n-1}=0,
\end{array}
\right.
\end{eqnarray*}
where $\bar{c}_{ij},1\le i,j\le n-1$, are positive
constants. Thus, we obtain by following the same argument for the
$(n-1)$ case that
$$
\bar{f}_i(x_1^*,\dots,x_{n-1}^*)=0,\ \ i=1,\dots,n-1.
$$
Therefore,
$$
f_i(x^*)= 0,\ \ \forall i\in \{1,2,\dots,n\}.
$$
\vskip 0.2cm
{\bf Case (c).}  Suppose that $\prod_{i=1}^{n}f_i(x^*)\neq 0$. We will
show that there is a contradiction. By symmetry, we need only
consider four different subcases as follows.
\vskip 0.2cm
{\bf Case (c1).} Suppose that
\begin{eqnarray*}
&&f_1(x^*)>0,\ f_i(x^*)<0,\ i=2,3,\dots,n.
\end{eqnarray*}
Similar to the case that $n=3$, we can find positive numbers
$\delta_1,\delta_2,\dots,\delta_{n-1}$ such that
\begin{eqnarray*}
&&x_1^*-\delta_1>0,x_2^*-\delta_2>0,\dots,x_{n-1}^*-\delta_{n-1}>0,\\
&&\frac{\delta_i}{\delta_1}=\frac{x_i^*}{x_1^*},\ i=2,3,\dots,n-1,
\end{eqnarray*}
and
\begin{eqnarray*}
&&f_1(x^*)>f_1(x_1^*-\delta_1,x_2^*-\delta_2,\dots,x_{n-1}^*-\delta_{n-1},x_{n}^*)>0,\\
&&f_2(x^*)<f_2(x_1^*-\delta_1,x_2^*-\delta_2,\dots,x_{n-1}^*-\delta_{n-1},x_{n}^*)<0,\\
&&\cdots\cdots\\
&&f_{n}(x^*)<f_{n}(x_1^*-\delta_1,x_2^*-\delta_2,\dots,x_{n-1}^*-\delta_{n-1},x_{n}^*)<0.
\end{eqnarray*}
It follows that
$$
F(x_1^*-\delta_1,x_2^*-\delta_2,\dots,x_{n-1}^*-\delta_{n-1},x_{n}^*)<F(x^*),
$$
which contradicts that $F$ reaches its minimum
at $x^*$.
\vskip 0.2cm
{\bf Case (c2).} Suppose that for $2\le i\le n-1$,
\begin{eqnarray*}
&&f_1(x^*)>0,f_2(x^*)>0,\dots,f_i(x^*)>0,\\
&&f_{i+1}(x^*)<0,\dots,f_{n}(x^*)<0.
\end{eqnarray*}
We fix $x_1^*,\dots,x_{i-1}^*$ and $x_{n}^*$. Similar to the
case that $n=3$, we can find positive numbers
$\delta_i,\cdots,\delta_{n-1}$ such that
\begin{eqnarray*}
x_i^*-\delta_i>0,\dots,x_{n-1}^*-\delta_{n-1}>0;\
\frac{\delta_j}{\delta_i}=\frac{x_j^*}{x_i^*},\  j=i+1,\dots,n-1,
\end{eqnarray*}
and
\begin{eqnarray*}
&&f_1(x^*)>f_1(x_1^*,\dots,x_{i-1}^*,x_i^*-\delta_i,\dots,x_{n-1}^*-\delta_{n-1},x_{n}^*)>0,\\
&&\cdots\cdots\\
&&f_i(x^*)>f_i(x_1^*,\dots,x_{i-1}^*,x_i^*-\delta_i,\dots,x_{n-1}^*-\delta_{n-1},x_{n}^*)>0,\\
&&f_{i+1}(x^*)<f_{i+1}(x_1^*,\dots,x_{i-1}^*,x_i^*-\delta_i,\dots,x_{n-1}^*-\delta_{n-1},x_{n}^*)<0,\\
&&\cdots\cdots\\
&&f_{n}(x^*)<f_{n}(x_1^*,\dots,x_{i-1}^*,x_i^*-\delta_i,\dots,x_{n-1}^*-\delta_{n-1},x_{n}^*)<0.
\end{eqnarray*}
It follows that
$$
F(x_1^*,\dots,x_{i-1}^*,x_i^*-\delta_i,\dots,x_{n-1}^*-\delta_{n-1},x_{n}^*)<F(x^*),
$$
which contradicts that $F$ reaches its minimum
at $x^*$.
\vskip 0.2cm
{\bf Case (c3).} Suppose that
\begin{eqnarray*}
&&f_1(x^*)<0,\ f_i(x^*)>0,\ i=2,3,\dots,n.
\end{eqnarray*}
Similar to the case that $n=3$, we can find positive numbers
$\delta_1,\delta_2,\dots,\delta_{n-1}$ such that
\begin{eqnarray*}
\frac{\delta_i}{\delta_1}=\frac{x_i^*}{x_1^*},\ i=2,3,\dots,n-1,
\end{eqnarray*}
and
\begin{eqnarray*}
&&f_1(x^*)<f_1(x_1^*+\delta_1,x_2^*+\delta_2,\dots,x_{n-1}^*+\delta_{n-1},x_{n}^*)<0,\\
&&f_2(x^*)>f_2(x_1^*+\delta_1,x_2^*+\delta_2,\dots,x_{n-1}^*+\delta_{n-1},x_{n}^*)>0,\\
&&\cdots\cdots\\
&&f_{n}(x^*)>f_{n}(x_1^*+\delta_1,x_2^*+\delta_2,\dots,x_{n-1}^*+\delta_{n-1},x_{n}^*)>0.
\end{eqnarray*}
It follows that
$$
F(x_1^*+\delta_1,x_2^*+\delta_2,\dots,x_{n-1}^*+\delta_{n-1},x_{n}^*)<F(x^*),
$$
which contradicts that $F$ reaches its minimum
at $x^*$.
\vskip 0.2cm
{\bf Case (c4).} Suppose that for $2\le i\le n-1$,
\begin{eqnarray*}
&&f_1(x^*)<0,f_2(x^*)<0,\dots,f_i(x^*)<0,\\
&&f_{i+1}(x^*)>0,\dots,f_{n}(x^*)>0.
\end{eqnarray*}
We fix $x_1^*,\dots,x_{i-1}^*$ and $x_{n}^*$. Similar to the
case that $n=3$, we can find positive numbers
$\delta_i,\dots,\delta_{n-1}$ such that
\begin{eqnarray*}
\frac{\delta_j}{\delta_i}=\frac{x_j^*}{x_i^*},\  j=i+1,\dots,n-1,
\end{eqnarray*}
and
\begin{eqnarray*}
&&f_1(x^*)<f_1(x_1^*,\dots,x_{i-1}^*,x_i^*+\delta_i,\dots,x_{n-1}^*+\delta_{n-1},x_{n}^*)<0,\\
&&\cdots\cdots\\
&&f_i(x^*)<f_i(x_1^*,\dots,x_{i-1}^*,x_i^*+\delta_i,\dots,x_{n-1}^*+\delta_{n-1},x_{n}^*)<0,\\
&&f_{i+1}(x^*)>f_{i+1}(x_1^*,\dots,x_{i-1}^*,x_i^*+\delta_i,\dots,x_{n-1}^*+\delta_{n-1},x_{n}^*)>0,\\
&&\cdots\cdots\\
&&f_{n}(x^*)>f_{n}(x_1^*,\dots,x_{i-1}^*,x_i^*+\delta_i,\dots,x_{n-1}^*+\delta_{n-1},x_{n}^*)>0.
\end{eqnarray*}
It follows that
$$
F(x_1^*,\dots,x_{i-1}^*,x_i^*+\delta_i,\dots,x_{n-1}^*+\delta_{n-1},x_{n}^*)<F(x^*),
$$
which contradicts that $F$ reaches its minimum
at $x^*$.\hfill\fbox

\vskip 0.5cm
\noindent {\bf Proof of Theorem \ref{thm2}}.
\vskip 0.2cm

{\it Case $n=2$.}
\vskip 0.2cm
Note that now equations (\ref{hsz1}) become
$$
q_{12}\frac{\alpha(2)}{\alpha(1)}=q_{21}\frac{\alpha(1)}{\alpha(2)}.
$$
Hence we can obtain a solution to (\ref{hsz1}) by defining $\frac{\alpha(2)}{\alpha(1)}=\sqrt{\frac{q_{21}}{q_{12}}}$.
\vskip 0.2cm
{\it Case} $n=3$.
\vskip 0.2cm
Equations (\ref{hsz1}) are equivalent to
\begin{eqnarray}\label{1}
\left\{
\begin{array}{l}
q_{12}\frac{\alpha(2)}{\alpha(1)}\beta(2)+q_{13}\frac{\alpha(3)}{\alpha(1)}\beta(3)=
q_{21}\frac{\alpha(1)}{\alpha(2)}\beta(2)+q_{31}\frac{\alpha(1)}{\alpha(3)}\beta(3),\\
q_{21}\frac{\alpha(1)}{\alpha(2)}\beta(1)+q_{23}\frac{\alpha(3)}{\alpha(2)}\beta(3)=
q_{12}\frac{\alpha(2)}{\alpha(1)}\beta(1)+q_{32}\frac{\alpha(2)}{\alpha(3)}\beta(3),\\
q_{31}\frac{\alpha(1)}{\alpha(3)}\beta(1)+q_{32}\frac{\alpha(2)}{\alpha(3)}\beta(2)=
q_{13}\frac{\alpha(3)}{\alpha(1)}\beta(1)+q_{23}\frac{\alpha(3)}{\alpha(2)}\beta(2).
\end{array}
\right.
\end{eqnarray}
Multiplying the first two equations by $\beta(1)/\beta(3)$ and $\beta(2)/\beta(3)$, respectively, and then adding them up,  we obtain the third equation. Define
$$
\frac{\alpha(2)}{\alpha(1)}=x_1,\ \frac{\alpha(3)}{\alpha(1)}=x_2.
$$
Thus, the first two equations of (\ref{1}) become
\begin{eqnarray}\label{2}
\left\{
\begin{array}{l}
a_{11}x_1+a_{12}x_2=b_{11}\frac{1}{x_1}+b_{12}\frac{1}{x_2},\\
a_{21}\frac{1}{x_1}+a_{22}\frac{x_2}{x_1}=b_{21}x_1+b_{22}\frac{x_1}{x_2},
\end{array}
\right.
\end{eqnarray}
where $a_{ij},b_{ij},1\le i,j\le 2$, are positive constants.

We define three continuous functions $f_1,f_2,F$ with the
domain $D_2:=\{x=(x_1,x_2)\in \mathbf{R}^2|x_1>0,x_2>0\}$ by
\begin{eqnarray*}
&&f_1(x)=(a_{11}x_1+a_{12}x_2)-\left(b_{11}\frac{1}{x_1}+b_{12}\frac{1}{x_2}\right),\\
&&f_2(x)=\left(a_{21}\frac{1}{x_1}+a_{22}\frac{x_2}{x_1}\right)-\left(b_{21}x_1+b_{22}\frac{x_1}{x_2}\right),\\
&&F(x)=f_1^2(x)+f_2^2(x).
\end{eqnarray*}
It is easy to see that the function $F$ has a minimum value at some point $x^*=(x^*_1,x^*_2)\in D_2$. Then, we have
\begin{eqnarray*}
&&\frac{\partial F}{\partial x_1}(x^*)=2f_1(x^*)\left(a_{11}+\frac{b_{11}}{(x^*_1)^2}\right)-2f_2(x^*)
\left(b_{21}+\frac{b_{22}}{x^*_2}+\frac{a_{21}}{(x^*_1)^2}+\frac{a_{22}x^*_2}{(x^*_1)^2}\right)=0,\\
&&\frac{\partial F}{\partial x_2}(x^*)=2f_1(x^*)\left(a_{12}+\frac{b_{12}}{(x^*_2)^2}\right)+2f_2(x^*)\left(\frac{a_{22}}{x^*_1}
+\frac{b_{22}x^*_1}{(x^*_2)^2}\right)=0.
\end{eqnarray*}
Since all the coefficients of the above equations are positive, we must have
$$
f_1(x^*)=0,\ f_2(x^*)=0.
$$
Hence there exists a positive solution $(x_1,x_2)$ to (\ref{2}) and therefore there exist $\alpha(i)>0, 1\le i\le 3$, such that
(\ref{hsz1}) holds.
\vskip 0.2cm
{\it Case} $n\ge 4$.
\vskip 0.2cm
Note that the last equation of (\ref{hsz1}) is implied by the first $(n-1)$ equations. If we define $\alpha(i)/\alpha(1)=x_{i-1}$ for $i=2,3,\dots,n$, then equations (\ref{hsz1}) become equations of the type (\ref{sys}). Therefore, the proof is completed by Theorem \ref{thm1}.\hfill\fbox

\section{Proof of Theorem \ref{hsz3}}\setcounter{equation}{0}

Let $\phi>0$ be a function on $E$. We define
$$
L^{\phi}_t=\frac{\phi(X_t)}{\phi(X_0)}\exp\left(-\int_0^t\frac{Q\phi}{\phi}(X_s)ds\right)\cdot
1_{\{t<\zeta\}},\ \ t\ge 0.
$$
$(L^{\phi}_t)_{t\ge 0}$ is a supermartingale of $X$. The upper bound (\ref{lower}) can be proved by following the standard argument (see \cite{I}). In the following, we will focus on the proof of the
lower bound (\ref{new1}).

Define $${\cal M}_0:=\{\mu\in {\cal P}_1(E): \mu(\{i\})>0,\
\forall i\in E\}.$$ Let $G$ be an open subset of ${\cal P}_1(E)$.
Denote by $m$ the measure on $E$ satisfying
$$m(\{i\})=\frac{1}{n},\ \ i\in E.$$ If $\delta>0$ is small enough,
then $(1-\delta)\mu+\delta m\in G\cap {\cal M}_0$ for each $\mu\in
G$. From the definition of $I(\mu)$, we find that
$$
I((1-\delta)\mu+\delta m)\le (1-\delta)I(\mu)+\delta I(m).
$$
Hence $\limsup_{\delta\rightarrow 0}[I((1-\delta)\mu+\delta m)]\le
I(\mu)$. Since $\mu\in G$ is arbitrary, $\inf_{\mu\in G}I(\mu)\ge
\inf_{\mu\in G\cap {\cal M}_0}I(\mu)$ and thus $\inf_{\mu\in
G}I(\mu)= \inf_{\mu\in G\cap {\cal M}_0}I(\mu)$. Therefore, to
prove (\ref{new1}), we need only prove that
\begin{equation}\label{new2}
\liminf_{t\rightarrow\infty}\frac{1}{t}\log P_{i}(L_t\in G,\
t<\zeta)\ge -\inf_{\mu\in G\cap {\cal M}_0}I(\mu),\ \ \forall i\in
E.
\end{equation}

Let $f$ be a function  on $E$. We define
$$
P^{\phi}_tf(i)=E_i(L^{\phi}_tf(X_t))=E_i\left(\frac{f(X_t)\phi(X_t)}{\phi(X_0)}\exp\left(-\int_0^t\frac{Q\phi}{\phi}(X_s)ds\right)\cdot
1_{\{t<\zeta\}}\right),\ \ t\ge 0.
$$
The generator of the semigroup $(P^{\phi}_t)_{t>0}$ is given by
$$
L^{\phi}f=\frac{Q(f\phi)}{\phi}-\frac{Q\phi}{\phi}f.
$$
That is, for any $i\in E$, we have
$$
L^{\phi}f(i)=\frac{Q(f\phi)(i)}{\phi(i)}-\frac{Q\phi(i)}{\phi(i)}f(i)=\sum_{j\in
E}\frac{q_{ij}\phi(j)}{\phi(i)}f(j)-\frac{\sum_{j\in
E}q_{ij}\phi(j)}{\phi(i)}f(i).
$$
Then, the matrix associated with $L^{\phi}$, denoted by
$Q^{\phi}=(q^{\phi}_{ij})_{i,j\in E}$, is given by
\begin{equation}\label{new3}
q^{\phi}_{ij}=\frac{q_{ij}\phi(j)}{\phi(i)},\ i,j\in E,i\neq j;\ \ q^{\phi}_{ii}=-\sum_{j\not=
i}\frac{q_{ij}\phi(j)}{\phi(i)},\ i\in E.
\end{equation}

Denote by $X^{\phi}$ the Markov chain associated with
$L^{\phi}$. By (\ref{new3}) and the assumption that $q_{ij}>0,\
i,j\in E,i\neq j$, we find that $X^{\phi}$ is an
ergodic Markov chain. Hence $X^{\phi}$ has a unique invariant
distribution, which is denoted by $\nu_{\phi}$. Note that
$$
P_{i}(L_t\in G,\
t<\zeta)=E^{\phi}_i\left(\frac{\phi(X_0)}{\phi(X_t)}\exp\left(\int_0^t\frac{Q\phi}{\phi}(X_s)ds\right);L_t\in
G\right).
$$
By the ergodicity of $X^{\phi}$, we obtain that
\begin{equation}\label{new5}\liminf_{t\rightarrow\infty}\frac{1}{t}\log P_{i}(L_t\in G,\
t<\zeta)\ge \sup_{\phi>0,\nu_{\phi}\in
G}\int_E\frac{Q\phi}{\phi}d\nu_{\phi},\ \ \forall i\in E.
\end{equation}

We define
$$
\Pi:=\{\mu\in {\cal P}_1(E):\mu=\nu_{\phi}\ {\rm for\ some}\
\phi>0\}.
$$
If we can prove the following claim
\begin{equation}\label{new4}
{\cal M}_0=\Pi,
\end{equation}
then we obtain by (\ref{new5}) that
\begin{eqnarray*}
-\inf_{\mu\in G\cap {\cal M}_0}I(\mu)&=&-\inf_{\mu\in G\cap {\cal
M}_0}\left\{-\inf_{u>0}\int_E\frac{Qu}{u}d\mu\right\}\\
&=&\sup_{\mu\in G\cap {\cal
M}_0}\left\{\inf_{u>0}\int_E\frac{Qu}{u}d\mu\right\}\\
&=&\sup_{\phi>0,\nu_{\phi}\in
G}\left\{\inf_{u>0}\int_E\frac{Qu}{u}d\nu_\phi\right\}\\
&\le&\sup_{\phi>0,\nu_{\phi}\in
G}\int_E\frac{Q\phi}{\phi}d\nu_{\phi}\\
&\le&\liminf_{t\rightarrow\infty}\frac{1}{t}\log P_{i}(L_t\in G,\
t<\zeta),
\end{eqnarray*}
and thus (\ref{new2}) is proved.

In the following, we will prove claim (\ref{new4}). Let $\mu\in {\cal M}_0$. We write
$$
d\mu=hdm,
$$
where $h$ is a function on $E$ satisfying $h(i)>0$
for each $i\in E$. To show that $\mu\in \Pi$, it is sufficient to show that there exists a function $\phi>0$ such that
$$
\int_E(L^{\phi}f)hdm=0,\ \ \forall f.
$$

Note that
\begin{eqnarray*}
&&\int_E(L^{\phi}f)hdm=0,\ \ \forall f\\
&\Leftrightarrow&\sum_{i\in E}(L^{\phi}f(i))h(i)=0,\ \ \forall f\\
&\Leftrightarrow&\sum_{i\in E}\left[\sum_{j\in
E}\frac{q_{ij}\phi(j)}{\phi(i)}f(j)-\frac{\sum_{j\in
E}q_{ij}\phi(j)}{\phi(i)}f(i)\right]h(i)=0,\ \ \forall f\\
&\Leftrightarrow&\sum_{j\in E}\sum_{i\in
E}\frac{q_{ji}\phi(i)}{\phi(j)}f(i)h(j)=\sum_{i\in
E}\frac{\sum_{j\in
E}q_{ij}\phi(j)}{\phi(i)}f(i)h(i),\ \ \forall f\\
&\Leftrightarrow&\phi(i)\sum_{j\in
E}q_{ji}\frac{h(j)}{\phi(j)}=\frac{h(i)}{\phi(i)}\sum_{j\in
E}q_{ij}\phi(j),\ \ \forall i\in E\\
&\Leftrightarrow&\phi^2(i)\sum_{j\in
E}q_{ji}\frac{h(j)}{\phi(j)}=h(i)\sum_{j\in
E}q_{ij}\phi(j),\ \ \forall i\in E.
\end{eqnarray*}
Hence, to show that $\mu\in \Pi$, it is sufficient to show that there exist
$\phi(i)>0$, $1\le i\le n$,  such that
\begin{equation}\label{trans}
h(i)\sum_{j\in E}q_{ij}\phi(j)=\phi^2(i)\sum_{j\in
E}q_{ji}\frac{h(j)}{\phi(j)},\ \ \forall i\in E.
\end{equation}
For $i\in E$, we define
$$
\beta(i)=\sqrt{h(i)},\ \
\alpha(i)=\frac{\phi(i)}{\beta(i)}.
$$
Then, equations (\ref{trans}) become equations (\ref{hsz1}). Since the existence of solutions to equations (\ref{hsz1}) is guaranteed by Theorem \ref{thm2}, the proof is complete. \hfill\fbox

\bigskip

{ \noindent {\bf\large Acknowledgments} \vskip 0.1cm  \noindent This work was supported by National Natural Science Foundation of China (Grant No. 11371191), Natural Sciences and Engineering Research Council of Canada (Grant No. 311945-2013), Natural Science Foundation of Hainan Province (Grant No. 117096), and Scientific Research Foundation for Doctors of Hainan Normal University.

\end{document}